%% file: RLDS-v3-arXiv.tex
\documentclass[12pt,draftcls,onecolumn]{IEEEtran}
\usepackage{fullpage,graphicx,psfrag,url,amsmath,bbm,enumerate,listings,amsfonts,amsthm,amssymb,cite}
\input defs.tex

\newcommand{\Rcal}{\mathcal{R}}
\newcommand{\Tcal}{\mathcal{T}}
\newcommand{\Ucal}{\mathcal{U}}
\newcommand{\Kcal}{\mathcal{K}}
\newcommand{\Dbold}{\mathbf{D}}
\newcommand{\Ebold}{\mathbf{E}}
\newcommand{\Ncal}{\mathcal{N}}
\def\Qt{\tilde{Q}}
\def\Vt{\tilde{V}}
\def\Dh{\hat{D}}
\def\e{\epsilon}
\def\s{\sigma}
\def\D{\Delta}

\nomenclature{$R$}{Number of markets for dispatch}
\nomenclature{$r$}{Index for a dispatch stage}
\nomenclature{$\Rcal$}{Index set for all dispatch stages}
\nomenclature{$s_r$}{Decision at $r$-th market}
\nomenclature{$x_r$}{Accumulated energy at $r$-th market}
\nomenclature{$x$}{Energy supply during the storage operation stages}
\nomenclature{$c_r$}{Price at $r$-th market}
\nomenclature{$c$}{Per-unit VOLL cost at each storage operation stage}
\nomenclature{$T$}{Number of storage operation stages}
\nomenclature{$t$}{Index for a storage operation stage}
\nomenclature{$\Tcal$}{Index set for storage operation stages}

\nomenclature{$L_t$}{Load at stage $t$}
\nomenclature{$W_t$}{Wind energy supply at stage $t$}
\nomenclature{$D_t$}{Deficit at stage $t$, \ie, $L_t - W_t$}
\nomenclature{$b_t$}{Storage level at stage $t$}
\nomenclature{$u_t$}{Storage operation at stage $t$}
\nomenclature{$\Ucal(b)$}{Feasible set for storage control given current storage level $b$}
\nomenclature{$F(b,u)$}{Storage system dynamic function}
\nomenclature{$Y_i$}{Information available at $i$-th stage, $i\in \Rcal \cup \Tcal$}
\nomenclature{$\phi_i$}{Control policy for stage $i$, $i\in \Rcal \cup \Tcal$}
\nomenclature{$B$}{Capacity of the storage device}
\nomenclature{$\lambda$}{Storage efficiency}
\nomenclature{$\mu$}{Recharging efficiency}
\nomenclature{$\nu$}{Discharging efficiency}
\nomenclature{$g(\cdot)$}{Terminal cost function for dispatch, \ie, the total cost of the storage operation stages}
\nomenclature{$J_r(x_r)$}{Cost-to-go for stage $r$}
\nomenclature{$\psi_r$}{Optimal dispatch threshold for stage $r$}
\nomenclature{$D_t^k$}{Effective deficit at stage $t$}
\nomenclature{$K_t$}{Width of the algebraic recombinant lattice at stage $t$}
\nomenclature{$\Kcal_t$}{Set of nodes of the algebraic recombinant lattice at stage $t$}
\nomenclature{$p_t^k$}{The probability of visiting $k$th node in the algebraic recombinant lattice at stage $t$}

\title{Risk Limiting Dispatch with Fast Ramping Storage}
\author{Junjie~Qin$^{1}$,
        Han-I~Su$^{2}$,
        and Ram~Rajagopal$^{3}$
\thanks{$^{1}$J. Qin is with the Institute for Computational and Mathematical Engineering and the Stanford Sustainable Systems Lab, Department of Civil and Environmental Engineering, Stanford University, CA, 94305. Email: jqin@stanford.edu.}%
\thanks{$^{2}$H. Su is with the Department of Electrical Engineering, Stanford University, CA, 94305. Email: hanisu@stanford.edu.}%
\thanks{$^{3}$R. Rajagopal is with the  Stanford Sustainable Systems Lab, Department of Civil and Environmental Engineering, Stanford University, CA, 94305. R. Rajagopal is supported by the Powell Foundation. Email: ramr@stanford.edu.}
\thanks{$*$ The first and second authors contributed equally to the paper. }
}

\begin{document}
\maketitle

\begin{abstract}
Risk Limiting Dispatch (RLD) was proposed recently as a mechanism that utilizes information and market recourse to reduce reserve capacity requirements, emissions and  achieve other system operator objectives. It induces a set of simple dispatch rules that can be easily embedded into the existing dispatch systems to provide computationally efficient and reliable decisions. Storage is emerging as an alternative to mitigate the uncertainty in the grid. This paper extends the RLD framework to incorporate fast-ramping storage. It developed a closed form threshold rule for the optimal stochastic dispatch incorporating a sequence of markets and real-time information. An efficient algorithm to evaluate the thresholds is developed based on analysis of the optimal storage operation. Simple approximations that rely on continuous-time approximations of the solution for the discrete time control problem are also studied. The benefits of storage with respect to prediction quality and storage capacity are examined, and the overall effect on dispatch is quantified. Numerical experiments illustrate the proposed procedures.
\end{abstract}

\section{Introduction}

The increased penetration of renewable generation in the grid increases the uncertainty that needs to be managed by the system operator \cite{NRELWest2010}. Existing  system dispatch procedures are designed assuming mild uncertainty, and utilize a worst-case schedule for generation by solving deterministic controls where random demands are replaced by their forecasted values added to `3$\sigma$'   \cite{Holttinen2008}, where $\sigma$ is the standard deviation of forecast error. Such rules require excess reserve capacity and result in increased emissions if $\sigma$ is large \cite{Hetzer2008, RBDEWV2012}.

An approach to mitigating the impact of increased uncertainty is to utilize energy storage \cite{Denholm2010}. Energy storage can be broadly classified into two groups depending on their scheduling characteristics \cite{Hadjipaschalis2009, Divya2009}: fast storage and slow storage. Fast storage can be utilized to mitigate intra-hourly variability of renewable supply. Slow storage can be utilized to transfer energy between different hours, so excess production can match peak demands.  In this paper we address the integration of fast storage into power system dispatch.

Existing approaches to stochastic dispatch rely on stochastic programming based on techniques like Monte Carlo scenario sampling that are hard to implement in large scale or do not incorporate market recourse \cite{Papavasiliou2011,Saric2009,Bouffard2008B,Morales2009,Conejo2010,IXJ2011}. Moreover, the optimal decisions can be difficult to interpret in the context of system operations. Recent work has proposed utilizing robust optimization formulations with uncertainty sets  \cite{Street2011, Jiang2012}, but they do not capture multiple recourse opportunities and can result in conservative dispatch decisions. Incorporating storage into these models results in additional complexity and decisions which are hard to analyze. Risk Limiting Dispatch (RLD) \cite{RBVW2011} was proposed as an alternative  to capture multiple operating goals and provide reliable and interpretable dispatch controls that can be readily incorporated in existing dispatch software.  RLD incorporates real-time forecast statistics and  recourse opportunities  enabling the evaluation of improvements in forecasting and market design \cite{RBDEWV2012}.  In this paper we develop RLD to incorporate fast storage.

Fast storage integration with renewables has been studied in a variety of scenarios. \cite{2011arXiv1109.3841S} examines the benefits of storage for renewable integration in a single bus model. Optimal contracting for co-located storage was studied in  \cite{Bitar2011,XGEE2012}, and the role of distributed storage was studied in  \cite{2011arXiv1110.4441K, GT2011}.  Recent independent work \cite{Harsha2012} addresses system operator (SO) dispatch of storage to mitigate net loads (scheduled load minus wind) to obtain analytic controls and expressions for the value of storage, but focusing only in the real-time dispatch.    We consider a SO dispatch process that utilizes a market with multiple recourse opportunities to evaluate the value of improved forecasts and impact of storage.  An easy-to-compute numerical dispatch algorithm is developed utilizing optimal control structural results and explicit formulae.  Analytic relationships between key quantities of interest are derived based on a continuous-time approximation of the storage problem.

The remainder of the paper is organized as follows.  Section~\ref{sec:problemStatement} states RLD with storage problem.
Section~\ref{sec:dispControl} establishes structural control results for dispatch.
The optimal storage operation and evaluation of the dispatch thresholds are studies in Section~\ref{sec:dtStorDispatch}.
An approximation scheme is derived in Section~\ref{sec:ctStorDispatch}. Numerical results are presented in Section~\ref{sec:numericalResults}. Section~\ref{sec:conclusion} concludes the paper.

\section{Problem Statement}\label{sec:problemStatement}
Grid operation is constrained so that supply must equal demand at each time instant. The system operator (SO) schedules conventional generation in a sequence of markets ahead of delivery time to ensure this constraint is met. Load and renewables are random, and only revealed at the delivery time interval. The goal of the operator is to find the optimal schedule and operation of grid resources given information about load and renewable generation. Fast storage can be used to smooth unpredicted variations of the net load (load minus renewables) in real-time.

\begin{figure}[h!]
\centering
\includegraphics[scale=0.55]{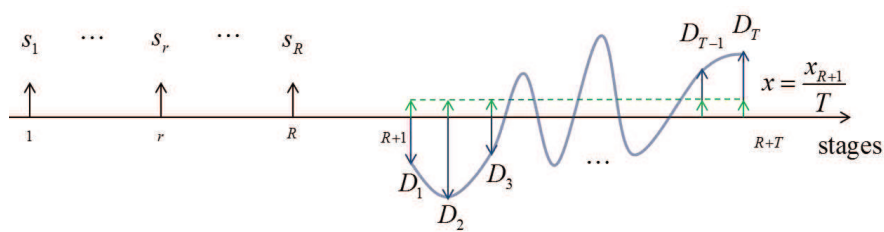}
\caption{Risk limiting dispatch with storage: The purchased conventional generation is supplied over the delivery time interval uniformly, while the wind generation and demand may vary over the time interval. A storage device is operated during the delivery time interval to minimize the terminal cost (penalty).}
\label{fig:RLDwStor}
\end{figure}

\subsection{Model Formulation}

There are $R$ markets, modeled as stages, ahead of the delivery time interval. For example, stage $1$ may occur $24$ hours ahead of the delivery time; stage $R$ occurs $15$ minutes ahead of the delivery time; other stages occur in between. At each stage, an operator makes a decision to purchase $s_r$ units of energy at stage $r$ for the delivery time interval. Note $s_r$ is a forward contract and may be interpreted as a contract for reserve capacity.  The price $c_r$ per unit of $s_r$ is known in advance. These reserve capacities are different in two respects: the $r$-th capacity must be available in shorter time than the $(r-1)$-th capacity, and their prices are different. The operator may also sell reserve capacity at various stages. In such case, $s_r < 0$.  We denote the index set for all dispatch stages as $\Rcal= \{1,2,\dots, R\}$, and use $r$ to denote one of the stages in $\Rcal$. Note the restriction of whether buying or selling is allowed at each stage is given ahead of time.

Further, to avoid trivial solutions, some constraints are imposed on the prices. For two buying dispatch stages $r_1,\ r_2\in \Rcal$ and $r_1 < r_2$, we require $0 < c_{r_1} <  c_{r_2}$,  i.e. price of purchasing power increases as the delivery deadline approaches. If $c_{r_1} > c_{r_2}$, it is worthwhile to defer the purchasing decision since more information is available at stage $r_2$ when the price is lower. Similarly, for two selling dispatch stages $r_1,\ r_2 \in \Rcal$ and $r_1 < r_2$, we require $c_{r_1} >  c_{r_2}$. Finally, to avoid arbitrage, for each buying stage $r_1$ and selling stage $r_2$, such that $r_1 < r_2$, we require $c_{r_1} > c_{r_2}$.

At each dispatch stage $r\in \Rcal$, three events occur. First information $Y_r$ is observed. In addition to the state at stage $r$, \ie, $x_r \in Y_r$, the information set could also contain signals that help the prediction of the demand and wind generation. Examples include weather forecast and sensor measurement data that are available at the time of stage $r$. Notice $Y_{r}\subset Y_{r+1}$, for all $r\in \Rcal$. Next a dispatch decision is made: The operator decides to purchase ($s_r \ge 0$) or sell ($s_r \le 0$) from the $r$-th market.  Lastly the total amount of power accumulated so far is computed
\begin{equation}
    x_{r+1} = x_r + s_r,\,\, r \in \Rcal.
\end{equation}
The energy accumulated in $R$ markets is supplied during a delivery time interval which is discretized into $T$ stages. Let  $\Tcal := \{R+1, R+2,\dots, R+T\}$ and use $t$ to denote each element in $\Tcal$. For each stage in the delivery time interval, the amount of energy supply from conventional generation is
$x= {x_{R+1}}/{T}$.

A random wind generation $W_t$ and a random load $L_t$ are realized for each $t\in \Tcal$. Let $D_t := L_t- W_t$ denote the net deficit at stage $t$. The deficit may be positive or negative. A monetary penalty is assessed to compensate the positive deficit or unmet demand. Different forms of the penalty will be discussed later in the section.
The information set is extended to stage $R+T$, \ie, for all $i \in \{2,\dots, R, R+1, \dots R+T\}$ we have $Y_{i-1}\subset Y_{i}$. We also define $Y_{R+T+1}$ as the information available at the end of the entire period and $Y_{R+T}\subset Y_{R+T+1}$.

In each stage $t\in \Tcal$, the storage operator can recharge $[u_t]_+$ or discharge $[u_t]_-$ units of energy subject to physical constraints of the storage device, where $[u_t]_+ = \max(u_t,0)$, $[u_t]_- = \min(u_t,0)$, and $u_t = [u_t]_+ + [u_t]_-$ is the variable representing the operation for storage at stage $t$. We denote the amount of energy stored in the storage device at stage $t$ as $b_t$, the transition function for the storage device as $F(b_t, u_t)$, and the feasible set for the discharging/recharging operations as $\mathcal{U}(b_t)$. We denote the terminal cost as $g(\mathbf{D},\mathbf{u},x)$, which will be specified in Subsection~\ref{sec:cost}, where $\mathbf{D} = (D_{R+1}, D_{R+2}, \dots, D_{R+T})$ and $\mathbf{u} = (u_{R+1}, u_{R+2}, \dots, u_{R+T})$. A control policy $\mathbf{\phi} = (\phi_1,\phi_2,\dots,\phi_{R+T})$ is a sequence of functions each of which maps the information available at current stage to the action at the same stage, \ie, $\phi_i$ is a $Y_i$-adapted function for every $i\in \Rcal\cup \Tcal$. We use $\mathbf{\phi}^D := (\phi_1,\phi_2,\dots,\phi_{R})$ to represent the dispatch policy and $\mathbf{\phi}^S := (\phi_{R+1},\phi_{R+2},\dots,\phi_{R+T})$ to represent the storage operation policy.

The RLD with storage problem can then be summarized as
\begin{subequations}\label{eq:RLDwStorProb}
\begin{align}
\mbox{minimize} \quad& \expec\left[\sum_{r\in \Rcal} c_rs_r +g(\mathbf{D},\mathbf{u},x)\right] \label{RLDwStorProbObej} \\
\mbox{subject to} \quad& s_r = \phi_r(\ycal_r), \label{RLDwStorProbC1}\\
                  & x_{r+1} = x_r + s_r, \;\; r\in \Rcal; \label{RLDwStorProbC2}\\
                  & u_t = \phi_t(\ycal_t), \label{RLDwStorProbC3}\\
                  & u_t \in \mathcal{U}(b_t), \label{RLDwStorProbC4}\\
                  & b_{t+1} = F(b_t, u_t), \;\; t\in \Tcal .\label{RLDwStorProbC5}
\end{align}
\end{subequations}

\subsection{Storage}
The storage device has  finite energy \emph{capacity} $B$.  Energy loss between stages is modeled using  the \emph{storage efficiency} $\lambda\in [0,1]$, so the contribution of $b_t$ units of stored energy at time $t+1$ is only $\lambda b_{t}$. Energy conversion loss is modeled with  parameters $\mu \in [0,1]$ and $\nu \in [0,1]$, denoting the \emph{recharging efficiency} and \emph{discharging efficiency}, respectively. If $u$ units of electric energy is recharged into the storage device, only $u'=\mu u$ units of energy will be actually stored due to energy conversion loss. Similarly, if $u'$ units of electric energy is discharged from the storage device, only $u=\nu u'$ unit of the energy can be used to meet the net deficit realized.

Typically, storage models also  consider ramping-rate constraints on charging and discharging \cite{QSVR2012}. Fast response grid level storage is rapidly becoming available with power to energy ratios of 40 to 50 kW per kWh utilizing Advanced Lithium-Ion blocks . A full charge or discharge of 1 kWh can be obtained in about 1.2 to 1.5 minutes of response time. If the dispatch discretization interval considered is larger than this,  the ramping constraints will not be active during operation.
Although some of our results can be generalized to cases with charge constraints (\eg, \cite{FRD2012,Harsha2012,Bitar2011}), focusing on a simpler model reveals deeper insight about the solution.  In ongoing work, we are devising a model for slow storage, which needs to account for the existence of multiple markets (each with multiple storage operation stages) with different timing constraints.

The dynamics of the storage model is captured by the transition function
\[
    F(b_t, u_t) = \lambda\left(b_t + \mu[u_t]_+ +\frac{1}{\nu} [u_t]_-\right),
\]
and the feasible action set is
\[
    \mathcal{U}(b_t) = \left\{u_t \middle| 0 \le [u_t]_+ \le \frac{1}{\mu}(B-b_t), - \nu b_t \le [u_t]_- \le 0\right\}.
\]
The feasible set for the vector $\mathbf{u}$ is denoted as $\mathcal{U}(\mathbf{b})$. We assume that the storage starts empty, \ie, $b_{R+1} =0$, and any energy remained in the storage after stage $R+T$ will be discarded at no cost/benefit consistent with operating policies for fast storage.

\subsection{Terminal Costs}
\label{sec:cost}
Different terminal costs lead to different dispatch goals: \\

\noindent{\bf Value of loss load (VOLL)}.  Let
\begin{align*}
g_t(D_t,u_t,x) = c[D_t - x + u_t]_+
\end{align*}
measure the shortfall to meet $D_t$ when $x$ units of power is supplied and $-u_t$ units of energy are withdrawn from the storage at stage $t$, where $c$ is the cost of unit shortfall.
The total terminal cost over the delivery time interval then is
\begin{align*}
g(\mathbf{D},\mathbf{u},x) = c\sum_{t \in \Tcal}[D_t - x + u_t]_+.
\end{align*}
Notice $g_t(D_t,u_t,x) = c [D_t - x + u_t]_+$ is convex in $x$, for all values of $D_t$ and $u_t$. \\

\noindent{\bf Loss of load probability (LOLP)}.  LOLP is in general defined as the probability of allocated supply not meeting the random deficit at the delivery time. For our model with a finite-length delivery time interval and storage, one way to define LOLP is
\begin{align*}
\prod_{t \in \Tcal} \prob\left( D_t + u_t >   x \middle| \ycal_{t}\right) \le \alpha,
\end{align*}
Here we use another definition which induces simple dispatch rule:
\begin{align*}\label{eq:riskcstr}
\prob\left( D_t + u_t >   x \middle| \ycal_{t}\right) \le \alpha_t,\,\, \forall t \in \Tcal,
\end{align*}
where $\alpha_t$ is the allowed LOLP at stage $t\in \Tcal$ which could be related to the allowed LOLP for the entire delivery time interval with \eg\  $\alpha_t = \alpha^{1/T}$.

A direct setup that achieves this goal is to use extended function definitions:
\begin{align*}
g(\mathbf{D},\mathbf{u},x) = \left\{\begin{array}{ll} 0  & \textrm{if } \prob\left( D_t + u_t >   x \middle| \ycal_{t}\right) \le \alpha_t,\,\, \forall t \in \Tcal,\\ \infty &  \textrm{otherwise. }\end{array}\right.
\end{align*}
Notice since the set
\[
    \left\{(\mathbf{u},x)\middle| \prob\left( D_t + u_t >   x \middle| \ycal_{t}\right) \le \alpha_t,\,\, \forall t \in \Tcal\right\}
\]
is convex, $g(\mathbf{D},\mathbf{u},x)$ is convex in $(\mathbf{u},x)$.\\

\noindent{\bf Frequency drop charge}.  In some scenarios, it is desirable to charge for the frequency deviations caused by unmet demand or excessive generation. A common assumption valid for small deviations of net demand is that the frequency deviation is linearly related to unmet demand, i.e. $\Delta f  = \alpha (D-x_{R+1})$.  If the deviation costs $\$ c$ per MW, then for a stage $t$
\begin{align*}
g_t(D_t,u_t,x)  = c\alpha |D_t-x + u_t|.
\end{align*}
In this case $g$ is also convex, but is not non-increasing. For all stages in the delivery time interval, we have
\begin{align*}
g(\mathbf{D},\mathbf{u},x) = c\alpha\sum_{t \in \Tcal}\left([D_t-x+u_t]_+ -  [D_t-x + u_t]_-\right).
\end{align*}

\section{Structural Result for Dispatch Control}\label{sec:dispControl}
The structure of the optimal dispatch depends on determining the optimal generation scheduling assuming that storage is operated optimally given a generator schedule. Based on the later assumption, a standard dynamic programming result reveals:
\begin{lemma}\label{thm:ctg}
The cost-to-go function for a dispatch stage $r\in \Rcal\cup \{R+1\}$ is
\begin{equation*}
    J_{R+1}(x_{R+1}) = \inf_{\mathbf{u}\in \mathcal{U}(\mathbf{b})} \expec \left\{g(\mathbf{D},\mathbf{u},x)\middle| Y_{R+1} \right\},
\end{equation*}
\begin{equation}\label{eq:ctg-general}
    J_r (x_r) = \inf_{s_r\in \mathcal{S}_r} \expec \left\{ c_r s_r + J_{r+1} (x_{r+1}) \middle| Y_r \right\} ,\,\, r\in \Rcal,
\end{equation}
where $\mathcal{S}_r = \{s| s\ge 0 \}$ if $r $ is a buying stage,  $\mathcal{S}_r = \{s| s\le 0 \}$ if $r$ is a selling stage.
Further, if $s_r^\star = \phi^\star_r(Y_r)$ minimizes the right hand side of~\eqref{eq:ctg-general} for each $x_r$ and $r$, then the dispatch policy $\mathbf{\phi}^{D\star}= (\phi^\star_1,\phi^\star_2,\dots,\phi^\star_{R}) $ is optimal.
\end{lemma}
Lemma \ref{thm:ctg} states that the policy minimizing cost-to-go functions is optimal. Note the per-stage terminal cost function $g_t(D_t,u_t,x)$ is jointly-convex in both $u_t$ and $x$. It follows that cost-to-go function for each stage $r\in\Rcal\cup \{R+1\}$ is convex:

\begin{proposition}\label{thm:convex}
The cost-to-go function $J_{r}(x_{r})$, for all $r\in \Rcal\cup \{R+1\}$ is convex in $x_{r}$ given $g(\mathbf{D},\mathbf{u},x)$ convex in $\mathbf{u}$ and $x$.
\end{proposition}

\begin{remark}
Proposition~\ref{thm:convex} is not restricted to the case of constant prices for dispatch. In fact, the convexity of the cost-to-go extends to the case where the price can depend on the dispatch, \ie, $c_r = c_r(s_r)$, as long as $c_r(\cdot)$ is a convex function for all $r\in \Rcal$.
\end{remark}

Based on these observations and principles from inventory control, the structure of the optimal control can be computed. This structure depends on the gradient of the cost function. When  cost functions are not differentiable (with respect to $x_r$), we use the notion of constrained subgradient denoted as $\grad$ in the sequel. Relying on Lemma~\ref{thm:ctg} and Proposition~\ref{thm:convex}, we are ready to give the main result of this section:
\begin{theorem}\label{thm:dispControlGen}
For each dispatch stage $r\in \Rcal$, the optimal dispatch is
\begin{equation}
    s_r^\star(x_r) = \begin{cases}
    [\psi_r - x_r]_+ &\mbox{ if $r$ is a buying stage,}\\
    [\psi_r - x_r]_- &\mbox{ if $r$ is a selling stage,}
    \end{cases}\label{eq:threshold}
\end{equation}
where  $\psi_r \in Y_r$ is a state independent variable that satisfies
\[
    c_r +\grad \hat{J}_{r+1}(\psi_r)=0,
\]
with $\hat{J}_{r+1}(x) = \expec [J_{r+1}(x)|Y_r]$. Thus $\psi_r$ is uniquely defined as $\psi_r = \grad \hat{J}_{r+1}^{-1}(-c_r)$.
\end{theorem}

Theorem~\ref{thm:dispControlGen} shows  the optimal dispatch is characterized by a sequence of thresholds $\psi_r$ for $r\in \Rcal$. This important practical feature of RLD \cite{RBV2011} generalizes to the case of storage with convex terminal costs. However, in RLD without storage thresholds can be precomputed given the probability structure of the net demand $D$ conditional on the information set $Y_r$. In the present case this is not possible in general as the net demand follows a stochastic process $D_t$, $t\in \Tcal$  instead of a single random variable representing the total net demand in the period. Moreover, the computation of the constrained subgradient of the terminal cost function coupled with minimization over feasible storage operations may not be analytically tractable.

\section{Storage Operation and Threshold Computation}\label{sec:dtStorDispatch}

The threshold structure derived in the Section~\ref{sec:dispControl} is valid for any choice of convex terminal cost function, but the actual threshold computation depends on the particular cost choice.  We focus on the VOLL cost in the reminder of the paper, but the analysis can be generalized to other costs in section~\ref{sec:cost}. For the VOLL cost,
the optimal control problem is
\begin{align*}
\mbox{minimize}\quad& \expec\left[\sum_{r\in \Rcal}c_rs_r + c\sum_{t \in \Tcal}[D_t - x + u_t]_+\right]  \\
\mbox{subject to} \quad& \eqref{RLDwStorProbC1}, \eqref{RLDwStorProbC2}, \eqref{RLDwStorProbC3}, \eqref{RLDwStorProbC4}, \eqref{RLDwStorProbC5}.
\end{align*}
Section~\ref{sec:dispControl}  solved this problem assuming storage is operated optimally. The remainder of the section derives a more explicit optimal control rule for storage  under the VOLL cost.  Based on it, an efficient algorithm for the constrained subgradient of the terminal cost-to-go function is developed, which simplifies the computation of the dispatch thresholds significantly.

\subsection{Storage Control}\label{sec:storControl}
The optimal storage operation problem, given $x_{R+1}$ units of energy accumulated in $R$ markets, is
\begin{align*}
\text{minimize} \quad  & \expec \left[c\sum_{t\in \Tcal} [D_t - x + u_t]_+ \right]  \\
\text{subject to} \quad & \eqref{RLDwStorProbC3}, \eqref{RLDwStorProbC4}, \eqref{RLDwStorProbC5}.
\end{align*}
The storage operation subproblem is again solved with dynamic programming.
\begin{lemma}[Optimal storage operation]\label{thm:storOper}
The terminal cost-to-go function is
\[
    J_{R+T+1}(x, b_{R+T}) = 0,
\]
and the cost-to-go for a storage operation stage $t\in \Tcal$ is
\[
    J_{t}(x,b_t) = \inf_{u_t \in \mathcal{U}(b_t)} \expec \left\{c[D_t - x + u_t]_+ + J_{t+1}(x, b_{t+1})\middle| Y_t \right\}.
\]
The cost-to-go function for each $t\in \Tcal\cup \{R+T+1\}$ is convex in $(x,b_t)$.
For $t\in \Tcal$, the optimal control policy for the storage operation is
\[
    u^\star_t(b_t) = \min\left\{[x-D_t]_+, \frac{1}{\mu} (B-b_t)\right\} - \min\left\{[D_t-x]_+, {\nu} b_t\right\},
\]
or equivalently in terms of recharging and discharging
\[
    [u^\star_t(b_t)]_+ = \min\left\{[x-D_t]_+, \frac{1}{\mu} (B-b_t)\right\},
    \quad
    [u^\star_t(b_t)]_- = - \min\left\{[D_t-x]_+, {\nu} b_t\right\}.
\]
\end{lemma}

\subsection{Threshold Computation}
The threshold can be computed combining analytic and algorithm approaches. Without loss of generality, we focus on the case of ideal storage ( $\lambda = \mu = \nu =1$) to simplify the notation. First a simple consequence of Lemma~\ref{thm:storOper} gives a  recursive formula for the expected total cost over all the storage operation stages (\ie, the terminal cost-to-go for dispatch after generators are scheduled):
\begin{corollary}[Terminal cost-to-go]\label{thm:terminalCostToGo}
The expected cost-to-go function at stage $R+1$ given the information at a dispatch stage $r$ is
\begin{align}\label{eq:JLongExpression} 
\hat{J}_{R+1}(x_{R+1}) = c \expec \left[\sum_{t\in \Tcal} [D_t-x-b_t^\star]\middle| Y_r \right],
\end{align}
where $x= x_{R+1}/T$ and the optimal storage level $b_t^\star$ is defined recursively as
\begin{equation}\label{eq:recursivebt}
b_{t+1}^\star = F_t(b_t^\star,u_t^\star) = B\wedge [x-D_t+b_t^\star]_+,
\end{equation}
with $b_{R+1}^\star = b_{R+1}=0$ and $a\wedge b = \min(a,b)$, for $t\in \Tcal$.
\end{corollary}
Eq.~\eqref{eq:JLongExpression} can be combined with prior results in \cite{RBV2011} to obtain the cost-to-go of other dispatch stages due to linearity of the cost structure. Also notice \eqref{eq:JLongExpression} gives a Monte Carlo based algorithm to evaluate the dispatch thresholds. In the reminder of this section, a closer investigation of the terminal cost-to-go and its subgradient is presented to devise more efficient algorithm for the threshold evaluation.

\begin{figure}[htbp]
\centering
\begin{tikzpicture}
\tikzstyle{level 1}=[sibling distance=24mm]
\tikzstyle{level 2}=[sibling distance=8mm]
\usepgflibrary{shapes.multipart}
\node [circle, draw,minimum size=6mm] {}
    child {node [circle, draw,minimum size=6mm]{}
        child {node [circle, draw,minimum size=6mm]{}}
        child {node [circle split,draw,circle split part fill={white,black!50},minimum size=6mm]{ \nodepart{lower} } }
        child {node [circle, draw, fill=black!50,minimum size=6mm]{}}}
    child {node [circle split,draw,circle split part fill={white,black!50},minimum size=6mm]{ \nodepart{lower} }
        child {node [circle, draw,minimum size=6mm]{}}
        child {node [circle split,draw,circle split part fill={white,black!50},minimum size=6mm]{ \nodepart{lower} } }
        child {node [circle, draw, fill=black!50,minimum size=6mm]{}}}
    child {node [circle, draw, fill=black!50,minimum size=6mm]{}
        child {node [circle, draw,minimum size=6mm]{}}
        child {node [circle split,draw,circle split part fill={white,black!50},minimum size=6mm]{ \nodepart{lower} } }
        child {node [circle, draw, fill=black!50,minimum size=6mm]{}}};
\end{tikzpicture}
\caption{Scenario tree for storage operation: $T=3$ example. The nodes depict the state of the storage device. A filled, half-filled, and unfilled node represents the case where the storage is full, between full and empty, and empty after the optimal control at the corresponding stage, respectively. The tree grows exponentially.}
\label{fig:storTree}
\end{figure}
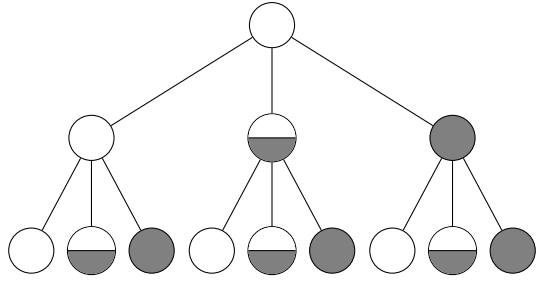

For this purpose, we first classify the states of the storage into three cases:  Empty ($b_t = 0$), full ($b_t = B$) and strictly in between ($0<b_t<B$). Then the proposed method works by calculating the probability of all possible sequences of states of the storage device. As an illustration, Fig.~\ref{fig:storTree} depicts the tree of possible storage states for the case $T=3$. Levels of the tree correspond to the storage operation stages, and nodes of the tree represent the state of storage device. Note the probability of visiting each node in the tree can be easily computed analytically.
However, this does not lead to a practical algorithm due to the curse of dimensionality. Because the number of nodes in each level of the tree grows exponentially as $t$ increases. In order to obtain a polynomial time algorithm, we introduce \emph{algebraic recombinant lattice}, an algebraic analog to
recombinant (or recombining) lattice, a technique widely used in finance applications \cite{Mun2002} and introduced to power systems and control community by \eg\ \cite{Papa2011}. In recombinant lattice model, the lattice (\ie, discretized state space) of dynamic programming has combined lattice points whenever two lattice points represents numeric values that are close enough, so that the growth of the state space is linear.
In a similar spririt, we combine the lattice points based on the algebraic forms in \eqref{eq:recursivebt}, or equivalently the algebraic form of \emph{effective deficit}, which is defined to be the difference between the realized deficit and the storage level for a particular case (node in Fig.~\ref{fig:storTree}).
\begin{proposition}[State space decomposition]\label{thm:stateSpaceDecompGen}
At stage $t\in \Tcal$, there are $K_t=2(t-R)-1$ algebraic forms of effective deficit $D_t^k$, which are defined recursively as
\begin{equation*}
    D_{t}^k = \begin{cases}
                    D_t & \mbox{ if }k = 1, \\
                    D_{t} - (x-D_{t-1}^{k-1}) &\mbox{ if } k \in \Kcal_{t}\backslash \{1,K_t\}, \\
                    D_{t} - B & \mbox{ if }k = K_t,
                    \end{cases}
\end{equation*}
where $k\in \Kcal_t =\{1,\dots, K_t\}$, or equivalently
\[
    D_{t}^k = \begin{cases}
    \sum_{i = 0 }^{k-1} D_{t-i} - (k-1)x & \mbox{ if }1 \le k \le t-R,\\
    \sum_{i = 0}^{K_t -k} D_{t-i} - (K_t - k)x - B & \mbox{ if }t-R < k \le K_t.
    \end{cases}
\]
The indicator of the event for each particular case to happen is
\begin{equation*}
    p_t^k = \begin{cases}
                \sum_{l \in \Kcal_{t-1}} p_{t-1}^l \ind\{D_{t-1}^l >  x\} & \mbox{ if }k =1, \\[.1 in]
                p_{t-1}^{k-1} \ind\{x-B < D_{t-1}^{k-1} \le  x\}                    & \mbox{ if }k \in \Kcal_{t} \backslash \{1, K_t\},\\[.1 in]
                \sum_{l \in \Kcal_{t-1}} p_{t-1}^l \ind\{D_{t-1}^l \le  x-B\} & \mbox{ if }k = K_t,
            \end{cases}
\end{equation*}
with $p_{R+1}^1 = 1$.
The probability of each of these events condition on information available at a dispatch stage $r$ is
$ \expec \left[p_t^k\middle| Y_r\right]. $
\end{proposition}

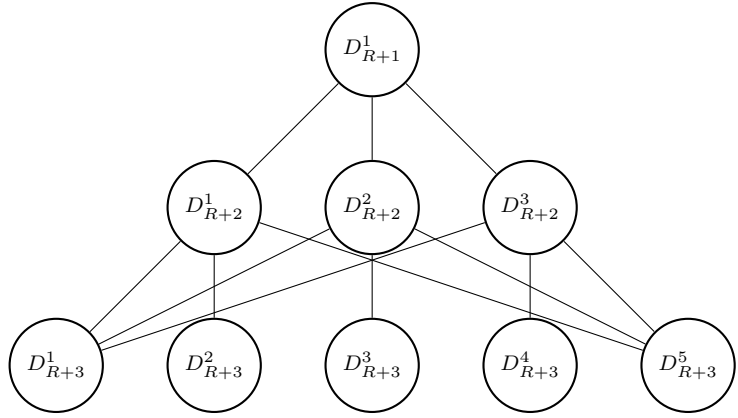
\begin{figure}[htbp]
\centering
\begin{tikzpicture}
  [scale=.7,auto=left,every node/.style={circle,shape=circle,draw,minimum size=6mm,thick}]
  \node (n11) at (1+3*2,1+3*2) {\scriptsize $D^1_{R+1}$};
  \node (n21) at (1+3*1,1+3)  {\scriptsize $D^1_{R+2}$};
  \node (n22) at (1+3*2,1+3)  {\scriptsize $D^2_{R+2}$};
  \node (n23) at (1+3*3,1+3) {\scriptsize $D^3_{R+2}$};
  \node (n31) at (1,1) {\scriptsize $D^1_{R+3}$};
  \node (n32) at (1+3*1,1)  {\scriptsize $D^2_{R+3}$};
  \node (n33) at (1+3*2,1)  {\scriptsize $D^3_{R+3}$};
  \node (n34) at (1+3*3,1) {\scriptsize $D^4_{R+3}$};
  \node (n35) at (1+3*4,1) {\scriptsize $D^5_{R+3}$};

  \foreach \from/\to in {n11/n21,n11/n22,n11/n23,n21/n31,n21/n32,n21/n35,n22/n31,n22/n33,n22/n35,n23/n31,n23/n34,n23/n35}
    \draw (\from) -- (\to);

\end{tikzpicture}
\caption{Evolution of effective deficit $D_t^k$: $T=3$ illustration. 
For each node $D_t^k$, its left and right children represents the case where starting with the effective deficit $D_t^k$, the storage ends up empty and full after optimal storage operation at stage $t+1$, respectively. The middle child represents the case that the storage level is strictly between $0$ and $B$.}
\label{fig:latticeDtk}
\end{figure}

The recombinant lattice based on the algebraic form of effective deficit $D_t^k$ is depicted in Fig.~\ref{fig:latticeDtk}. We will also refer the $k$-th node on the $t$-th level (corresponding to $D_t^k$ in the figure) node $(t,k)$. Notice only the nodes originally corresponding to cases where the storage is empty or full are combined since they share the same expression of $D_t^k$. Now the number of nodes in each level grows linearly as a function of $t$, as desired. Furthermore, the recursive definition of the indicator of visiting each node characterizes the condition for each case. The expected terminal cost-to-go and its subgradient with respect to $x$ can then be expressed in terms of the effective deficit $D_t^k$ and indicators $p_t^k$.

\begin{theorem}[Terminal cost-to-go and subgradient]\label{thm:gradientComp}
The expected cost-to-go function at stage $R+1$ given the information at a dispatch stage $r$ is
\begin{align}
\hat{J}_{R+1}(x_{R+1})= &c \expec \left[\sum_{t\in \Tcal} \sum_{k \in \Kcal_t}[D_t^k - x]_+p_t^k\middle|Y_r\right] \label{eq:JhatDtkForm}\\
=& c \expec \left[\sum_{t\in \Tcal} \sum_{k \in \Kcal_t}\sum_{j \in \Ncal_{t}^k}[D_t^k - x]\ind \{(\underline{\Dbold_t^k})_j < \Dbold_t \le  (\overline{\Dbold_t^k})_j \}\middle|Y_r\right] \nonumber
\end{align}
where
$    \Dbold_{t} = \left( D_{R+1}^1,  D_{R+2}^1,  D_{R+2}^2,  D_{R+2}^3,\dots, D_{t}^1, \dots,   D_{t}^{K_t}
    \right)^\mathrm{T}$,
$(\underline{\Dbold_t^k})_j$ and $(\overline{\Dbold_t^k})_j$ are columns of matrix $\underline{\Dbold_t^k}\in \real^{(t-R)^2 \times |\Ncal_t^k|}$ and $\overline{\Dbold_t^k}\in \real^{(t-R)^2 \times |\Ncal_t^k|}$, which are chosen such that
$    \ind\{D_t^k>  x\}p_t^k = \sum_{j \in \Ncal_{t}^k}\ind \{ (\underline{\Dbold_t^k})_j < \Dbold_{t} \le  (\overline{\Dbold_t^k})_j\}$,
with each element in $\Ncal_t^k$ representing a path from the root node to node $D_t^k$ in the lattice.

The constrained subgradient for the cost-to-go function at stage $R+1$ given the information at a dispatch stage $r$ is
\begin{align*}
   & \grad \hat{J}_{R+1}(x_{R+1}) \\
    =& -\frac{c}{T} \sum_{t\in \Tcal} \sum_{k \in \Kcal_t} \sum_{j \in \Ncal_{t}^k} \Bigg\{ \left[k \wedge(K_t+1-k)\right] \expec  \left[ \ind \{(\underline{\Dbold_t^k})_j < \Dbold_{t} \le  (\overline{\Dbold_t^k})_j \}\middle|Y_r \right] \\
    +& \sum_{i=1}^{(t-R)^2} \expec \Bigg[ (D_t^k - x)\left(\frac{d (\overline{\Dbold_t^k})_{i,j}}{dx} \ind \{(\Dbold_{t})_i=(\overline{\Dbold_t^k})_{i,j}\}-\frac{d (\underline{\Dbold_t^k})_{i,j}}{dx} \ind \{(\Dbold_{t})_i=(\underline{\Dbold_t^k})_{i,j}\}\right)\\
    &\quad \quad \quad \quad\ind \{(\underline{\Dbold_t^k})_{-i,j}< (\Dbold_{t})_{-i}\le (\overline{\Dbold_t^k})_{-i,j}\}\Bigg| Y_r \Bigg]\Bigg\},
    \label{eq:subgradientGeneral}
\end{align*}
where $(\Dbold_{t})_i$, $(\underline{\Dbold_t^k})_{i,j}$ and $(\overline{\Dbold_t^k})_{i,j}$ are $i$-th entries of column vector $\Dbold_{t}$, $(\underline{\Dbold_t^k})_{j}$ and $(\overline{\Dbold_t^k})_{j}$, respectively; $(\Dbold_{t})_{-i}$, $(\underline{\Dbold_t^k})_{-i,j}$ and $(\overline{\Dbold_t^k})_{-i,j}$ are the remaining parts of the corresponding vectors.
\end{theorem}

\subsection{Gaussian Dispatch}\label{sec:gaussianDispatch}
Theorem~\ref{thm:gradientComp} works for general deficit processes. In practice, information about the net load is given by forecasts of load and wind and the expected variance of these forecasts. Utilizing the predicted deficit,  the dispatch can be simplified by considering the forecast error random variable of the net load. The prediction errors of the deficit process can be assumed to be  Gaussian random variables as observed in various studies (\eg\, \cite{Harsha2012, RBV2011, Makarov2009}). In particular we consider the following form of the forecast
\begin{equation}\label{eq:forecastModel}
    D_t = \hat{D_t}(Y_r) + \epsilon_t(Y_r), \quad \forall t\in \Tcal,
\end{equation}
where $\epsilon_t(Y_r)\sim N\left(0,\sigma_t^2(Y_r)\right)$ is independently distributed for each $t$. For each dispatch stage $r\in \Rcal$, the forecast $\hat{D_t}(Y_r)$  and forecast error $\epsilon_t(Y_r)$ depend on the information $Y_r$.
A typical pattern of this dependence is that as the delivery time approaches, the variance in prediction error decreases due to the accumulated information that is collected by \eg\, wind speed sensors around the wind farm. This dependence is captured by inputting different variance values for the prediction error at different dispatch stages. Note for each fixed $t\in \Tcal$, $\epsilon_t(Y_r)$ may not be independent across dispatch stages (for different values of index $r$). Since the calculation in the reminder of this paper applies to each dispatch stage, we write $\epsilon_t$ and $\hat{D}_t$ directly when index $r$ is clear from the context, and omit the dependence on $Y_r$ to simplify the notation.

The independence of prediction errors simplifies the evaluation of the probability of visiting each node. We now give the updated version of Proposition~\ref{thm:stateSpaceDecompGen}.

\begin{proposition}[State space decomposition: Gaussian errors]\label{thm:stateSpaceDecompIndependent}
With forecast model~\eqref{eq:forecastModel}, the predicted effective deficit $\hat{D}_t^k$  is
\begin{equation}\label{eq:DtkIndependent}
    \hat{D}_{t}^k = \begin{cases}
    \sum_{i = 0 }^{k-1} \hat{D}_{t-i} - (k-1)x & \mbox{ if }1 \le k \le t-R,\\
    \sum_{i = 0}^{K_t -k} \hat{D}_{t-i} - (K_t - k)x - B & \mbox{ if }t-R < k \le K_t,
    \end{cases}
\end{equation}
and the prediction error in the effective deficit is
\begin{equation}\label{eq:etk}
    \epsilon_t^k =  \sum_{i = 0 }^{h} \epsilon_{t-i},
\end{equation}
where $h= (k-1)\wedge (K_t-k)$ is the depth of the node $(t,k)$ from the closest boundary nodes.

For the state corresponding to node $(t,k)$ in the algebraic recombinant lattice, denote the probability to visit the node as $p_t^k$, visit the node and move to its left child as $\overleftarrow{p_t^k}$, visit the node and move to its middle child as $\overline{p_t^k}$, and visit the node and move to its right child as $\overrightarrow{p_t^k}$. Following recursion holds for these quantities. Starting with $p_{R+1}^1=1$,
\begin{equation}
    p_t^k = \begin{cases}
                \sum_{l \in \Kcal_{t-1}} \overleftarrow{p_{t-1}^l}        & \mbox{ if }k =1, \\[.1 in]
                \overline{p_{t-1}^{k-1}}     & \mbox{ if }k \in \Kcal_{t} \backslash \{1, K_t\},\label{eq:ptkIndependent} \\[.1 in]
                \sum_{l \in \Kcal_{t-1}} \overrightarrow{p_{t-1}^l}       & \mbox{ if }k = K_t,
            \end{cases}
\end{equation}
\begin{equation}
    \overleftarrow{p_t^k}=p_{t-h}^{k-h}\prob \left(\underline{\Ebold_{t-1}^{k-1}}< \Ebold_{t-1}^{k-1}\le \overline{\Ebold_{t-1}^{k-1}}, \,\, \epsilon_t^k>   \overline{\epsilon_t^k}\right), \label{eq:ptkLeft}
\end{equation}
\begin{equation}
    \overline{p_t^k}=p_{t-h}^{k-h}\prob \left( \underline{\Ebold_{t}^{k}} < \Ebold_{t}^{k} \le \overline{\Ebold_{t}^{k}}\right),\label{eq:ptkMiddle}
\end{equation}
\begin{equation}
    \overrightarrow{p_t^k}=p_{t-h}^{k-h}\prob \left(  \underline{\Ebold_{t-1}^{k-1}}< \Ebold_{t-1}^{k-1}\le \overline{\Ebold_{t-1}^{k-1}} ,\,\, \epsilon_t^k\le \underline{\epsilon_t^k}\right),\label{eq:ptkRight}
\end{equation}
where
$\underline{\epsilon_t^k}=x-B-\hat{D}_t^k$,
$  \overline{\epsilon_t^k} = x-\hat{D}_t^k$,
$\Ebold_t^k = \left(\epsilon_{t-h}^{k-h}, \epsilon_{t-h+1}^{k-h+1},  \dots,  \epsilon_{t}^{k} \right)^\mathrm{T}$,
$\underline{\Ebold_t^k} = \left(\underline{\epsilon_{t-h}^{k-h}}, \underline{\epsilon_{t-h+1}^{k-h+1}},  \dots,  \underline{\epsilon_{t}^{k}} \right)^\mathrm{T}$,
and $\overline{\Ebold_t^k} = \left(\overline{\epsilon_{t-h}^{k-h}}, \overline{\epsilon_{t-h+1}^{k-h+1}},  \dots,  \overline{\epsilon_{t}^{k}} \right)^\mathrm{T}$.
\end{proposition}

Given that $\epsilon_t$ is independently distributed zero-mean Gaussian, $\epsilon_t^k$ is also zero-mean Gaussian, whose variance can be easily computed. It follows that $\Ebold_t^k$ is a zero-mean multivariate Gaussian and its distribution function can be evaluated provided the its covariance matrix which again is available from the definition of $\epsilon_t^k$. Proposition \ref{thm:stateSpaceDecompIndependent} allows us to evaluate the expected terminal cost-to-go and its subgradient explicitly.

\begin{lemma}[Terminal cost-to-go and subgradient: Gaussian errors]\label{thm:offline}
The expected terminal cost is
\begin{align*}
    \hat{J}_{R+1}(x_{R+1})
    = c \sum_{t\in \Tcal} \sum_{k \in \Kcal_t}   \left[\hat{D}_t^k + \mu\left(\Ebold_t^k;\begin{bmatrix}\underline{\Ebold_{t-1}^{k-1}}\\ \dotfill\\ \overline{\epsilon_t^k}\end{bmatrix} , \begin{bmatrix}\overline{\Ebold_{t-1}^{k-1}}\\ \dotfill\\ \infty\end{bmatrix}\right)_{h+1} -x\right] \overleftarrow{p_t^k} ,
\end{align*}
where $\mu\left(X; \underline{X}, \overline{X}\right)$ is the mean vector of the truncated Gaussian with mean and variance equal to that of $X$, and truncation interval  $\left[\underline{X}, \overline{X}\right)$. Here the second term in the bracket is the mean of the last entry of $\Ebold_t^k$ within the corresponding interval.

The expected constrained subgradient is
\begin{equation}
 \grad \hat{J}_{R+1}(x_{R+1})
    = -\frac{c}{T} \sum_{t\in \Tcal} \sum_{k \in \Kcal_t}   \left[k \wedge(K_t+1-k)\right] \overleftarrow{p_t^k}.
    \label{eq:subgradientIndependent}
\end{equation}
\end{lemma}

Relying only on the evaluation of Gaussian distribution function, Proposition~\ref{thm:stateSpaceDecompIndependent} and Lemma~\ref{thm:offline} give an analytical tractable approach to calculate the expected total cost for the delivery interval with storage operation. It provides an efficient approach to compute the dispatch thresholds for prediction model \eqref{eq:forecastModel}.
We also note this result is of interest in other applications of energy storage, where benefits of storage need to be quantified as cost-saving over a finite horizon.
We can in fact analyze several practical special cases of model~\eqref{eq:forecastModel} and point out cases where the thresholds can be computed off-line:
\begin{enumerate}
\item $\sigma_t \equiv \sigma$. In this case, the prediction errors are i.i.d. Gaussian random variables. This type of prediction models are typically favored by power engineers because there may not be enough data to estimate different variances for the prediction error at different storage operation stages. Further, this results in a simpler implementation for the calculation. Note in terms of the analytic derivation, this assumption does not lead to further simplification.
\item $\hat{D}_t \equiv \hat{D}$. In this case, the only forecast available is one nominal value for the deficit over the entire delivery time interval. This simplifies the form of $D_t^k$ in Equation~\eqref{eq:DtkIndependent}. From a practical perspective, this assumption allows the computation of the thresholds to be conducted off-line.
This simplifies the dispatch procedure tremendously, so that stochastic dispatch may be carried out in a similar fashion as conventional deterministic dispatch.
\item $\hat{D}_t \equiv \hat{D}$ and $\sigma_t \equiv \sigma$. This is the simplest model in which thresholds can be computed off-line. It also requires extremely few data to estimate model parameters. However, this model may be too simple to represent the fluctuation of the deficit process over the delivery time interval.
\end{enumerate}

\section{Approximate algorithm for dispatch}\label{sec:ctStorDispatch}
In this section, we consider the continuous-time operation of energy storage and propose an approximate algorithm for estimating dispatch thresholds. Before introducing the continuous-time model, we first reformulate the discrete-time counterpart. Without loss of generality, we assume $c = 1$. Let
\begin{align*}
  V_t & = \sum_{\tau = R+1}^t [D_\tau - x + u_\tau]_+, \\
  Q_t & = \sum_{\tau = R+1}^t [D_\tau - x + u_\tau]_-,
\end{align*}
denote the cumulative VOLL cost and cumulative curtailment up to time $t$, respectively. Suppose that $V_t$ and $Q_t$, $t \in \Tcal$, are adapted to information $Y_t$, $t \in \Tcal$. Then the charging and discharging operation of energy storage is uniquely specified by
\begin{align*}
  u_t & =
  \begin{cases}
    (V_t - V_{t-1}) + (Q_t - Q_{t-1}) - (D_t - x) & \text{if}\ t > R+1, \\
    V_{R+1} + Q_{R+1} - (D_{R+1} - x) & \text{if}\ t = R+1,
  \end{cases}
\end{align*}
and is also adapted to $Y_t$, $t \in \Tcal$. The stored energy can also be expressed in terms of $V_t$ and $Q_t$:
\begin{align*}
  b_{t+1}
  & = b_t + u_t \\
  & = b_{R+1} + \sum_{\tau = R+1}^t u_\tau \\
  & = b_{R+1} - \sum_{\tau = R+1}^t (D_t - x) + V_t + Q_t.
\end{align*}
Now we can reformulate the optimal storage operation problem with $V_t$ and $Q_t$ as control variables, that is,
\begin{align*}
  \text{minimize} \quad & \expec\left[  V_{R+T+1} \right] \\
  \text{subject to} \quad
  & b_{t+1} = b_{R+1} - \sum_{\tau=R+1}^t (D_\tau - x) + V_t + Q_t, \\
  & 0 \leq b_{t+1} \leq B, \\
  & V_t \geq V_{t-1} \geq \ldots \geq V_{R+1} \geq 0, \\
  & Q_t \leq Q_{t-1} \leq \ldots \leq Q_{R+1} \leq 0, \\
  & (V_t, Q_t) = \phi_t(Y_t).
\end{align*}
Although the feasible set allowing $(V_t - V_{t-1})(Q_t - Q_{t-1}) < 0$ is larger than the feasible set of the original problem, it is easy to see that the alternative control variables
\begin{align*}
  \Vt_t & = V_t - \min\{V_t - V_{t-1}, -(Q_t - Q_{t-1})\}, \\
  \Qt_t & = Q_t + \min\{V_t - V_{t-1}, -(Q_t - Q_{t-1})\},
\end{align*}
yield the same stored energy $b_{t+1}$ and lower cost. Thus, the reformulated optimization problem is equivalent to the original problem.
Under the optimal policy in Lemma~\ref{thm:storOper}, the cumulative VOLL cost $V_t$ increases only if storage is empty, that is, $b_{t+1} = 0$, and the cumulative curtailment $Q_t$ increases only if storage is full, that is, $b_{t+1} = B$.

With the above reformulation, we are ready to introduce the continuous-time model. Assume that the delivery time is a continuous-time interval $\Tcal_C := [R+1, R+T+1]$. Assume that given information set $Y_{R+1}$, the cumulative net deficit process $\mathbf{D}$ is a $(\Dh(Y_{R+1})/T, \s_{R+1}/\sqrt{T})$ Brownian motion, that is, $D_{t}$ is a Gaussian random variable with mean $\Dh(Y_{R+1})(t-R-1)/T$ and variance $\s_{R+1}^2 (t-R-1)/T$, where the adjustment in $t$ is due to the starting time.
The cumulative VOLL cost $V_t$ is adapted to information $Y_t$, continuous, and non-decreasing with $V_{R+1} = 0$. The cumulative curtailment $Q_t$ is adapted to information $Y_t$, continuous, and non-increasing with $Q_{R+1} = 0$.
Then the stored energy at time $t$ is equal to
\begin{align*}
  b_t & = b_{R+1} - [D_t - x (t-R-1)] + V_t + Q_t
\end{align*}
for $t \in \Tcal_C$.
Under the optimal policy, $V_t$ increases only if $b_t = 0$, and $Q_t$ decreases only if $b_t = B$. The stored energy process $b_t$ is a reflected Brownian motion. We will approximate the total VOLL cost by the product of the long-term average VOLL cost and the delivery interval length. To find the long-term average cost, we use the properties of reflected Brownian motion in the following Lemma.

\begin{lemma}[\cite{H1985}]
Let $Z_t$ be a $(\mu, \s)$ Brownian motion with $Z_0 = 0$ and
\begin{align*}
  b_t = Z_t + V_t + Q_t
\end{align*}
be a reflected Brownian motion in $[0, B]$ such that $V_t$ and $Q_t$ are adapted to the filtration induced by $Z_t$ and satisfy
\begin{enumerate}
  \item $V_t$ is continuous and non-decreasing with $V_0 = 0$,
  \item $V_t$ increases only when $b_t = 0$,
  \item $Q_t$ is continuous and non-increasing with $Q_0 = 0$, and
  \item $Q_t$ decreases only when $b_t = B$.
\end{enumerate}
The long term averages of $V_t$ and $Q_t$ are equal to
\begin{align*}
  \lim_{t \to \infty} \frac{1}{t} \expec[V_t]
  & =
  \begin{cases}
    \displaystyle \frac{\mu}{e^{2\mu B/\s^2} - 1} & \text{if}\ \mu \neq 0, \\
    \displaystyle \frac{\s^2}{2 B} & \text{if}\ \mu = 0,
  \end{cases} \\
  \lim_{t \to \infty} \frac{1}{t} \expec[Q_t]
  & =
  \begin{cases}
    \displaystyle - \frac{\mu}{1 - e^{-2\mu B/\s^2}} & \text{if}\ \mu \neq 0, \\
    \displaystyle - \frac{\s^2}{2 B} & \text{if}\ \mu = 0,
  \end{cases}
\end{align*}
respectively. The steady-state probability density function of $Z_t$ is equal to
\begin{align*}
  f_Z(z) & =
  \begin{cases}
    \displaystyle \frac{2\mu}{\s^2} \left( \frac{e^{2\mu z / \s^2}}{e^{2\mu B/\s^2} - 1} \right) & \text{if}\ \mu \neq 0, \\
    \displaystyle \frac{1}{B} & \text{if}\ \mu = 0,
  \end{cases}
\end{align*}
for $0 \leq z \leq B$.
\end{lemma}

Using the above lemma, we can approximate the VOLL cost for general $c$ and its first-order derivative by
\begin{subequations}\label{eq:meanPenalty}
\begin{align}
  J_{R+1}(x_{R+1})
  & = c \expec[V_T]
  \approx c T \lim_{t \to \infty} \frac{1}{t} \expec[V_t] \\
  & = \frac{c \s_{R+1}^2}{2 B} h\left( \frac{2 B}{\s_{R+1}^2} (x_{R+1} - \Dh(Y_{R+1})) \right), \nonumber\\
  \grad J_{R+1}(x_{R+1})
  & \approx c h'\left( \frac{2 B}{\s_{R+1}^2} (x_{R+1} - \Dh(Y_{R+1})) \right).
\end{align}
\end{subequations}
where
\begin{align*}
  h(x) & =
  \begin{cases}
    \displaystyle \frac{x}{e^x - 1} & \text{if}\ x \neq 0, \\
    \displaystyle 1 & \text{if}\ x = 0,
  \end{cases} \\
  h'(x) & =
  \begin{cases}
    \displaystyle \frac{(1-x)e^x - 1}{(e^x - 1)^2} & \text{if}\ x \neq 0, \\
    \displaystyle - \frac{1}{2} & \text{if}\ x = 0.
  \end{cases}
\end{align*}
\begin{remark}
Formulae~\eqref{eq:meanPenalty} reveal the role played by storage explicitly. An important observation is that scaling $B$ and $\sigma_{R+1}^2$ by the same constant does not affect $J_{R+1}(x_{R+1})$ and its derivative. That is, a system with more fluctuate wind (deeper penetration) and large storage can have the same terminal cost and thus dispatch thresholds with the another system with less fluctuate wind and small storage, given the ratio $B/\sigma_{R+1}^2$ is fixed. This quantifies the notion ``storage firms the wind'' in the context of dispatch.
\end{remark}

The approximate VOLL cost is convex. Thus, the approximate dispatch policy is still characterized by~\eqref{eq:threshold} except that the thresholds $\psi_r = \D_r + \Dh(Y_r)$ and $\D_r$ satisfies
\begin{align*}
    c_r & = c_{r+1} \left( 1 - \prob\{ \D_r > \D_{r+1} + \e_r \}  \right) \\
    & \qquad + c_{r+2} \left( \prob\{ \D_r > \D_{r+1} + \e_r \} - \prob\{ \D_r > \D_{r+1} + \e_r,\ \D_r > \D_{r+2} + \e_r^{r+1} \}  \right) + \ldots \\
    & \qquad + c_R ( \prob\{ \D_r > \D_{r+1} + \e_r,\ \ldots,\ \D_r > \D_{R-1} + \e_r^{R-2} \} \\
    & \qquad\qquad - \prob\{ \D_r > \D_{r+1} + \e_r,\ \ldots,\ \D_r > \D_R + \e_r^{R-1} \} ) \\
    & \qquad - c \expec\left[ \ind_{\{ \D_r > \D_{r+1} + \e_r,\ \ldots,\ \D_r > \D_R + \e_r^{R-1} \}} h'\left( \frac{2B}{\s_{R+1}^2} (\D_r - \e_r^R) \right)\right].
\end{align*}

\section{Numerical results}\label{sec:numericalResults}

\subsection{Setup}
We utilize the published forecast performance curve from Red Electrica Espana ( Fig.~\ref{fig:prelim}\subref{fig:err}) to compare the costs of various dispatch policies. Let $\s(t)$ be the standard deviation of the $t$-hours-ahead forecast error. The error $\e_r$ explained from stage $r-1$ to stage $r$ is assumed to be a Gaussian random variable with zero mean and variance $\s_r^2 = \s(t_{r-1})^2 - \s(t_r)^2$. We assume that at $t = 0.25$, the forecast error of the mean of the deficit contributes $20\%$ of the error variance, and thus $\s_{R+1}^2 = 0.8 \s(0.25)^2$. For the discrete-time model, the number of storage operation time intervals is $|\Tcal| = 60$.

We consider a 3-stage dispatch with day-ahead, hour-ahead, and 15-minutes-ahead stages. The prices of purchasing energy are suggested by published average energy prices in California. We set the day-ahead price to $\$52$ per MWh, the hour-ahead price to $\$60$ per MWh, the 15-minutes-ahead price to $\$72$ per MWh, and the VOLL to $\$1000$ per MWh. The mean of the deficit $D$ is normalized and is between $-1$ and $+1$. For a policy $\phi$, we will estimate the cost $J_\phi(D, B)$ by $2000$ Monte Carlo runs of forecast errors.

\subsection{Comparing dispatch approaches}
In addition to the optimal dispatch policy in Section~\ref{sec:dtStorDispatch} and the approximate algorithm in Section~\ref{sec:ctStorDispatch}, we also consider the following two dispatch approaches as benchmarks. The $3\s$-rule assumes $\D_r = 3\s(t_r)$. The ideal policy is the optimal dispatch given a perfect forecast and is given by the following linear program,
\begin{align*}
  \text{minimize} \qquad & c_1 x_1 + c \sum_{t \in \Tcal} [D_t - x + u_t]_+ \\
  \text{subject to} \qquad
  & b_{t+1} = b_t + u_t, \\
  & - b_t \leq u_t \leq B - b_t.
\end{align*}
Since a perfect forecast is available in this case, it is always optimal to make all purchase at the day-ahead market when the price is lowest. Denote the cost of the ideal policy by $J_0(D, B)$. For any policy $\phi$, $J_\phi(D, B) \geq J_0(D, B)$, and the difference is the {\it integration cost} of policy $\phi$: $C_I = J_\phi(D, B) - J_0(D, B)$ \cite{RBDEWV2012}.

\begin{figure}
\subfloat[Forecast error curve]{
\footnotesize
\psfrag{xlabel}[c][c]{Time horizon (hour)}
\psfrag{ylabel}[c][c]{Forecast error}
\psfrag{legenderr}[r][r]{Forecast error}
\psfrag{legenderrc}[r][r]{Dispatch stages}
\includegraphics{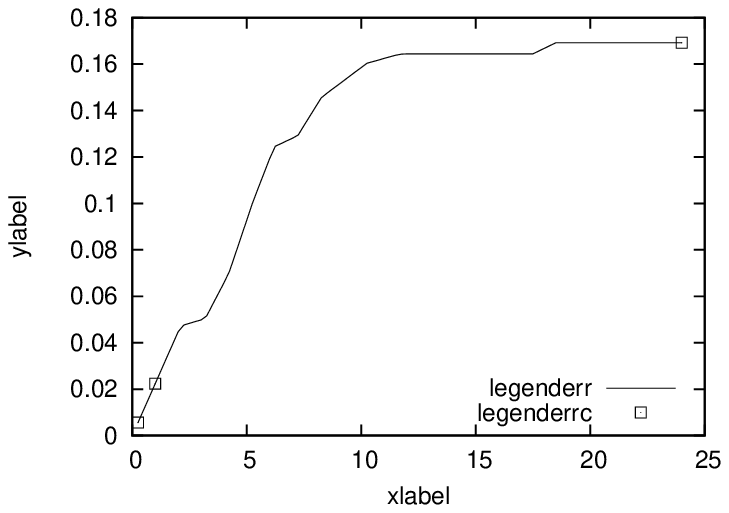}
\label{fig:err}
}
\subfloat[Storage operation costs]{
\footnotesize
\psfrag{xlabel}[c][c]{$x_{R+1}$}
\psfrag{ylabel}[c][c]{$J_{R+1}$}
\psfrag{legendjcs}[r][r]{Continuous-time, $B = 0.001$}
\psfrag{legendjds}[r][r]{Discrete-time, $B = 0.001$}
\psfrag{legendjcl}[r][r]{Continuous-time, $B = 0.01$}
\psfrag{legendjdl}[r][r]{Discrete-time, $B = 0.01$}
\includegraphics{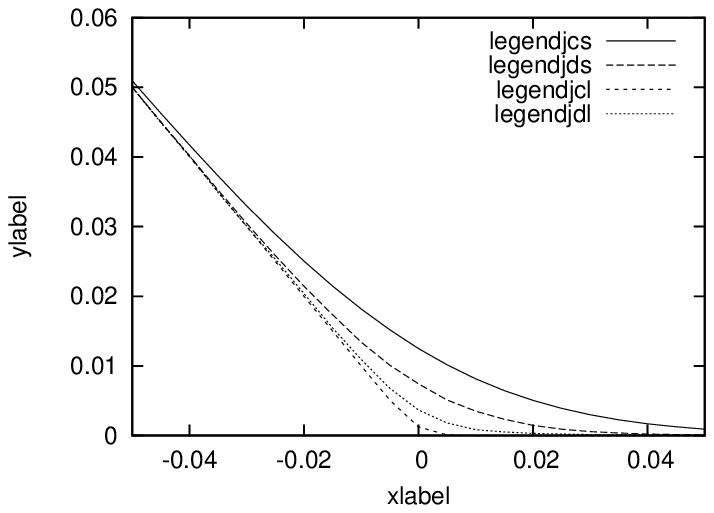}
\label{fig:cost}
}
\caption{Illustration of the forecast error curve from Red Electrica Espana and the storage operation costs for the discrete-time and continuous-time models.}
\label{fig:prelim}
\end{figure}

Fig.~\ref{fig:prelim}\subref{fig:cost} shows the storage operation costs for the discrete-time model and the approximate continuous-time model. The approximate model overestimates the storage operation cost for small storage capacity since the discrete-time model does not consider the cost caused by the variation within each time interval. For large storage capacity, the continuous-time model underestimates the storage operation cost since it assumes that the probability distribution of the initial stored energy is steady-state distribution instead of zero assumed in the discrete-time model.

\begin{figure}
\centering
\subfloat[Costs]{
\footnotesize
\psfrag{xlabel}[c][c]{$D$}
\psfrag{ylabel}[c][c]{$J_\phi(D, B)$}
\psfrag{legend3sigma}[r][r]{$3 \sigma$}
\psfrag{legendconti}[r][r]{Continuous-time}
\psfrag{legenddiscrete}[r][r]{Discrete-time}
\psfrag{legendideal}[r][r]{Ideal}
\includegraphics{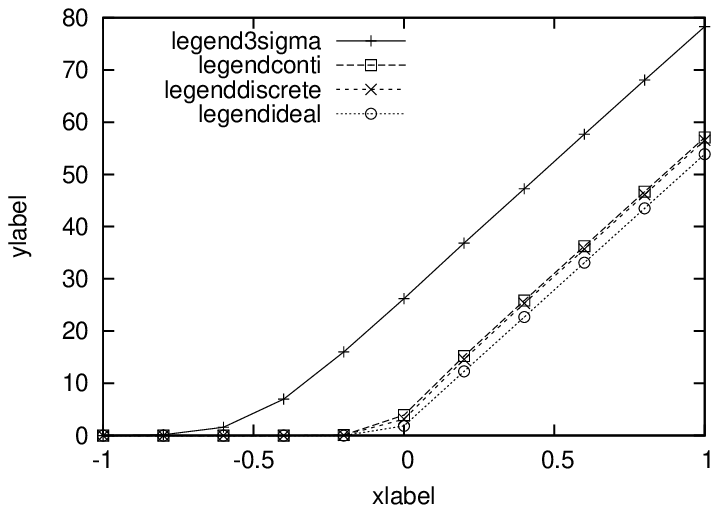}
\label{fig:costD}
}
\subfloat[Integration costs]{
\centering
\footnotesize
\psfrag{xlabel}[c][c]{$D$}
\psfrag{ylabel}[c][c]{$C_I$}
\psfrag{legend3sigma}[r][r]{$3 \sigma$}
\psfrag{legendconti}[r][r]{Continuous-time}
\psfrag{legenddiscrete}[r][r]{Discrete-time}
\includegraphics{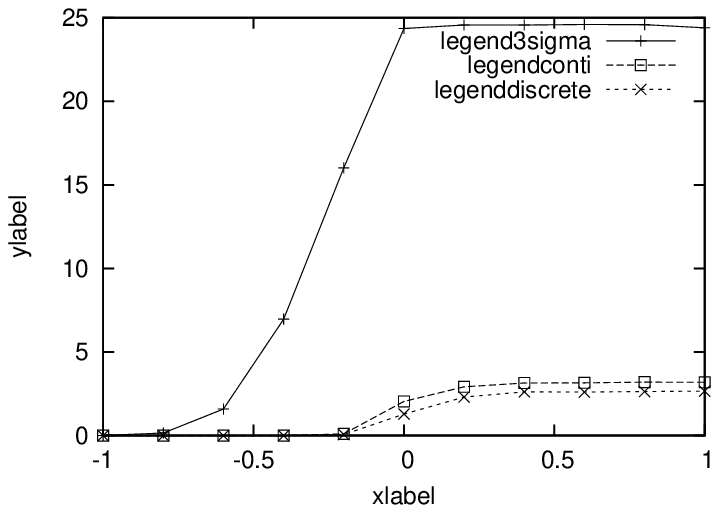}
\label{fig:costI}
}
\caption{The costs and integration costs for the $3 \sigma$ strategy, the optimal dispatch policy, the approximate policy, and the ideal policy for $B = 0.001$. }
\label{fig:costDI}
\end{figure}

In Fig.~\ref{fig:costDI}\subref{fig:costD}, we compare the cost $J_{\phi}(D, B)$ of the $3 \sigma$ strategy, the optimal dispatch policy, the approximate policy, and the ideal policy for $B = 0.001$. The $3 \sigma$ strategy has the highest cost. The cost of the approximate policy is slightly higher than the optimal cost. The integration costs with respect to the cost of the ideal policy are shown in Fig.~\ref{fig:costDI}\subref{fig:costI}.

\begin{figure}
\centering
\footnotesize
\psfrag{xlabel}[c][c]{$B$}
\psfrag{ylabel}[c][c]{$J_\phi(D, B)$}
\psfrag{legendconti}[r][r]{Continuous-time}
\psfrag{legenddiscrete}[r][r]{Discrete-time}
\includegraphics{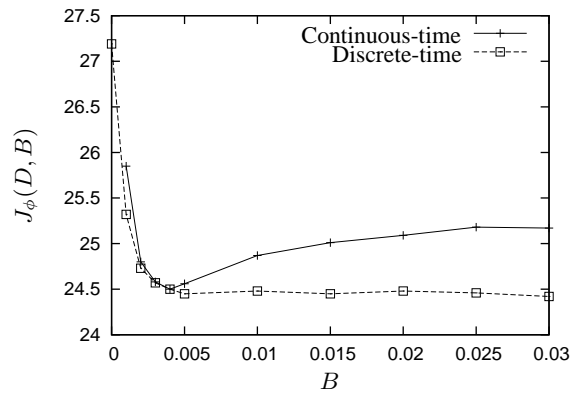}
\caption{The costs of the optimal policy and the approximate policy for $D = 0.4$.}
\label{fig:costB}
\end{figure}

Fig.~\ref{fig:costB} shows the cost $J_{\phi}(D, B)$ of the optimal dispatch policy and the approximate policy for $D = 0.4$. As we observe in Fig.~\ref{fig:cost}, the approximate model is not suitable for very small and very large storage capacities and thus has higher costs in those regimes.

\section{Conclusion}\label{sec:conclusion}
The paper extends Risk Limiting Dispatch to incorporate fast ramping storage. The structural properties of the optimal dispatch are studied, and demonstrate that  optimality is achieved by following a simple threshold rule. The optimal storage operation policy is given in closed form. Explicit formulae for evaluating the total expected cost over the delivery interval are obtained and efficient algorithms for computing the dispatch threshold using this cost estimates are also obtained.

A simpler continuous time approximation to storage operation results in a simple expression for  terminal cost-to-go as a function of the storage capacity $B$ and the deficit process variance $\sigma_{R+1}^2$. The relationship quantifies the notion that the storage smoothes the wind. The algorithms are illustrated and compared using numerical results. In ongoing work, we are extending the method to include slow storage, which requires modeling multiple simultaneous market decisions. We would also like to investigate incorporation of ramping constraints for the slow storage problem. Finally, the effects of network congestion in a scenario with storage can be considered.

\section*{Acknowledgement}
The authors would like to acknowledge the support from TomKat Center,  Powell Foundation Fellowship.  The authors would like to thank Pravin Varaiya, Eilyan Bitar and Kameshwar Poolla for invaluable discussions.

\appendices
\section{Proof for Section~\ref{sec:dispControl}}

\noindent{\bf Proof of Proposition~\ref{thm:convex}}\\
We first state and prove a standard result.
\begin{proposition}[\cite{HS2003}]\label{thm:HS2003}
Let $X$ be a nonempty set with $A_x$ a nonempty set for each $x\in X$. Let $C=\{(x,y):y \in A_x, x\in X\}$, let $J$ be a real-valued function on $C$ and define
\begin{align*}
f(x)=\inf \{J(x,y):y\in A_x\}\  , x\in X.
\end{align*}
If $C$ is a convex set and $J$ is a convex function on $C$, then $f$ is a convex function on any convex subset of $X^\star=\{x:x\in X,f(x)>-\infty\}$.
\end{proposition}
\begin{IEEEproof}
Pick $x_1$ and $x_2$ in $X^\star$ so $f(x_1)>-\infty$, $f(x_2)>-\infty$. Then for all $\gamma >0$ there are $y_1$ and $y_2$ with $(x_i,y_i)\in C$, $i=1$ and $2$, such that $f(x_i)+\gamma> J(x_i,y_i)$. Pick $t\in (0,1)$, and let $(x,y)=t (x_1,y_1)+(1-t)(x_2,y_2)$, which is in $C$ because $C$ is convex. Now
\begin{align*}
tf(x_1)+(1-t)f(x_2)\ge tJ(x_1,y_1)+(1-t)J(x_2,y_2)-\gamma\ge J(x,y)-\gamma \ge f(x)-\gamma,
\end{align*}
with the second inequality due to convexity of $J$ on $C$. Letting $\gamma \rightarrow 0$ yields the convexity of $f$.
\end{IEEEproof}

Given $g(\mathbf{D}, u, x)$ convex in $\mathbf{u}$ and $x$, we have $\expec \left\{ g(\mathbf{D}, u, x) \middle|Y_{R+1}\right\}$ convex in $\mathbf{u}$ and $x$. Since $\mathcal{U}(\mathbf{b})$ is an affine set (and therefore is convex), by Proposition~\ref{thm:HS2003}, we have $J_{R+1}(x_{R+1})$ convex in $x_{R+1}$.

Suppose $J_{r+1}(x_{r+1})$ is convex in $x_{r+1}$. We proceed to prove $J_{r}(x_{r})$ is convex in $x_r$. Note
\[
    J_r (x_r) = \inf_{s_r\in \mathcal{S}_r} \expec \left\{ c_r s_r + J_{r+1} (x_{r+1}) \middle| Y_r \right\}.
\]
Since
\[
    c_r s_r + J_{r+1} (x_{r}+s_{r}),
\]
is convex in $s_{r}$ and $x_{r}$, and the conditional expectation preserves the convexity, by invoking Proposition~\ref{thm:HS2003}, we have $J_{r}(x_{r})$ is convex in $x_{r}$.

Therefore by induction, we have $J_{r}(x_{r})$ convex in $x_{r}$ for all $r\in \Rcal\cup \{R+1\}$.
\bigskip

\noindent{\bf Proof of Theorem~\ref{thm:dispControlGen}}\\
By Proposition~\ref{thm:convex}, we have $J_{r}(x_r)$ is convex for $r\in \Rcal\backslash\{R+1\}$. Further $-c_r \in \grad \hat{J}_{r+1}(\psi_r)$ by the definition of constrained subgradient. Since $\hat{J}_{r+1}$ is $Y_r$-adapted, we have $\psi_r$ is $Y_r$-adapted.  Given
\[
    \hat{J}_{r+1}(x) - \hat{J}_{r+1}(\psi_r) \ge -c_r (x-\psi_r),
\]
we have
\[
    c_r x + \hat{J}_{r+1}(x) \ge c_r \psi_r+ \hat{J}_{r+1}(\psi_r).
\]
Therefore
\[
    s_r^\star = \psi_r - x_r
\]
if $s_r\in \mathcal{S}_r$. This relation gives the optimal threshold for stages $r\in \Rcal $. For purchase only or sell only stages, we show Equation~\eqref{eq:threshold} gives the optimal dispatch. Consider a purchase stage, if $\psi_r - x_r \ge 0$, the constraint is not tight and $s_r^\star = \psi_r - x_r$. Otherwise, we show $s_r^\star = 0$. Suppose $s_r^\star = \tilde{s}_r^\star \not = 0$, there exists $0 \le \alpha \le 1$ such that $\alpha \tilde{s}_r^\star + (1-\alpha) (\psi_r - x_r) =0$. However, given $\tilde{J}_r(x_r,s_r)= \expec\left[c_r s_r +J_{r+1}(x_r + s_r)\middle| Y_r\right] $ convex in $s_r$, we have
\[
    \tilde{J}_r(x_r,\tilde{s}_r^\star)< \tilde{J}_r(x_r,0) = \tilde{J}_r(x_r,\alpha \tilde{s}_r^\star + (1-\alpha) (\psi_r - x_r) ) \le \alpha \tilde{J}_r(x_r,\tilde{s}_r^\star) + (1-\alpha)\tilde{J}_r(x_r,\psi_r-x_r),
\]
where the first inequality is based on the assumption that $\tilde{s}_r^\star$ is a minimizer of $\tilde{J}_r(x_r,s_r)$ while $0$ is not, the last inequality is due to the convexity. Consequently
\[
    \tilde{J}_r(x_r,\tilde{s}_r^\star)<\tilde{J}_r(x_r,\psi_r-x_r),
\]
which is clearly a contradiction since $\psi_r-x_r$ is the minimizer to the unconstrained problem. Therefore $s_r^\star = [\psi_r - x_r]_+$. Similarly for the sell only stage, $s_r^\star = [\psi_r - x_r]_-$.

\section{Proof for Section~\ref{sec:dtStorDispatch}}
\noindent{\bf Proof of Lemma~\ref{thm:storOper}}\\
We need to prove the convexity of the cost-to-go function and the optimality of the proposed control rule.
\begin{itemize}
    \item Notice $J_{R+T+1}(x, b_{R+T})$ is convex in $(x, b_{R+T})$.

        Suppose $J_{t+1}(x,b_{t+1})$ is convex  in $(x, b_{t+1})$.  We proceed to prove $J_{t}(x,b_{t})$ is convex  in $(x,b_{t})$. Note
        \[
            J_{t}(x,b_t) = \inf_{u_t \in \mathcal{U}(b_t)} \expec \left\{c[D_t - x + u_t]_+ + J_{t+1}(x,b_{t+1}) \middle| Y_t \right\}.
        \]
        Since
        \[
            [ D_t - x + u_t ]_+ + J_{t+1} \left( x, \lambda \left( b_t + \mu [u_t]_+ + \frac{1}{\nu} [ u_t ]_- \right) \right)
        \]
        is convex in $(x,b_t, u_t)$ , and the conditional expectation preserves the convexity, by invoking Proposition~\ref{thm:HS2003}, we have $J_{t}(x, b_{t}) $ is convex in $(x,b_{t})$. 

        Therefore by induction, we have $J_{t}(x, b_{t})$ is convex in $(x,b_{t})$ for all $t\in \Tcal\cup \{R+T+1\}$.
    \item The optimal control policy of storage is a standard result. See {\bf Remark 4.3} in \cite{Bitar2011} for intuitional explanation, and \cite{2011arXiv1109.3841S} for detailed proof.
\end{itemize}
\bigskip

\noindent{\bf Proof of Corollary~\ref{thm:terminalCostToGo}}\\
We prove this result by induction over the length of the storage operation problem. For the base case, let $T=1$, we have
\[
    \hat{J}_{R+1}(x_{R+1};T=1) = c \expec \left\{ [D_{R+1}-x]_+ \middle| Y_r \right\},
\]
\[
    b^\star_{R+1} = 0, \mbox{ and } b^\star_{R+2} = B\wedge [x-D_t+b^\star_{R+1}].
\]
Suppose the expression for $\hat{J}_{R+1}(x_{R+1})$ and $b^\star_{t}$ holds for $T= l$, and consider the case $T=l+1$. By the optimal storage control rule, we have
\[
    \hat{J}_{R+1}(x_{R+1};T=l+1) = c \expec \left\{ \hat{J}_{R+1}(x_{R+1};T=l) + [D_{R+l+1} - x - b^\star_{R+l+1}]_+ \middle| Y_r\right\}.
\]
Invoking the induction hypothesis on the recursive formula for the sequence $b^\star_{t}$, for $T=l$, whose last term gives $b^\star_{R+l+1}$, we have
\[
    \hat{J}_{R+1}(x_{R+1};  T=l+1) = c \expec \Bigg\{ \hat{J}_{R+1}(x_{R+1};T=l) +
    \left[D_{R+l+1} - x - B\wedge \left[x-D_{R+l}+b^\star_{R+l} \right]_+\right]_+ \Bigg| Y_r\Bigg\}.
\]
Plugging in the expression of $ \hat{J}_{R+1}(x_{R+1};T=l)$ yields the desired result on $\hat{J}_{R+1}(x_{R+1})$, and the expression of $b^\star_{R+l+2}$, which is the only additional term in the sequence $b^\star_{t}$, for $T=l+1$, follows from the optimal storage control rule.
\bigskip

\noindent{\bf Proof of Proposition~\ref{thm:stateSpaceDecompGen} and Theorem~\ref{thm:gradientComp}}\\
For the sake of the limited space, we prove Proposition~\ref{thm:stateSpaceDecompGen} in the context of Theorem~\ref{thm:gradientComp}. The general proof of Proposition~\ref{thm:stateSpaceDecompGen} can be done similarly by induction.
Consequently there are two items to prove:
\begin{itemize}
\item Proposition~\ref{thm:stateSpaceDecompGen} holds, \ie, $\hat{J}_{R+1}(x_{R+1})$ in Corollary~\ref{thm:terminalCostToGo} can be expressed as
\begin{align*}
\hat{J}_{R+1}(x_{R+1})= c \expec \left[\sum_{t\in \Tcal} \sum_{k \in \Kcal_t}[D_t^k - x]_+p_t^k\middle|Y_r\right].
\end{align*}
\begin{IEEEproof}
Equivalent to the form in Corollary~\ref{thm:terminalCostToGo}, we have
\begin{align}
\hat{J}_{R+1}(x_{R+1}) = c \expec &\Bigg\{[D_{R+1}-x]_+ + \left[D_{R+2}-x-B \wedge \left[x-D_{R+1}\right]_+\right]_+  \nonumber\\
& + \left[D_{R+3}-x-B \wedge \left[x-D_{R+2}+B\wedge\left[x-D_{R+1}\right]_+\right]_+\right]_+ + \dots \label{eq:JhatContReveal} \\
&
+ \left[D_{R+T} - x - B\wedge \left[x-D_{R+T-1}+B\wedge \left[\dots\left[x-D_{R+1}\right]_+\dots\right]_+\right]_+\right]_+\Bigg| Y_r \Bigg\}. \nonumber
\end{align}
Denote the expected penalty that will occur at stage $t\in \Tcal$ as  $V_t$, \ie,
\[
    V_t = c \expec \Bigg\{\left[D_{t} - x - B\wedge \left[x-D_{t-1}+B\wedge \left[\dots\left[x-D_{R+1}\right]_+\dots\right]_+\right]_+\right]_+ \Bigg| Y_r \Bigg\}.
\]
We have
$
    \hat{J}_{R+1}(x_{R+1}) = \sum_{t\in \Tcal} V_t.
$
We then need to prove
$
    V_t = c \expec \Bigg\{\sum_{k\in \Kcal_t} [D_t^k -x]_+ p_t^k \Bigg| Y_r\Bigg\},
$
\ie,
\begin{align}
    & \expec \Bigg\{\sum_{k\in \Kcal_t} [D_t^k -x]_+ p_t^k \Bigg| Y_r\Bigg\} \label{eq:DtkEquivalence}\\
    &\qquad =\expec \Bigg\{\left[D_{t} - x - B\wedge \left[x-D_{t-1}+B\wedge \left[\dots\left[x-D_{R+1}\right]_+\dots\right]_+\right]_+\right]_+ \Bigg| Y_r \Bigg\} \nonumber
\end{align}
which will prove the expression for $\hat{J}_{R+1}(x_{R+1}) $ by summing up terms corresponding to each storage operation stages. By observation, the equation above holds if under the event indicated by $p_t^k$ we have
\[
    D_{t}- B\wedge \left[x-D_{t-1}+B\wedge \left[\dots\left[x-D_{R+1}\right]_+\dots\right]_+\right]_+ = D_t^k.
\]
We prove this statement by induction. The base case holds since $\Kcal_{R+1} = \{1\}$ , $p_{R+1}^1 = 1$ and $D_{R+1}^1 = D_{R+1}$. Suppose the result holds for stage $t-1$ and consider the stage $t$.
If the event indicated by $p_t^1$ holds, \ie, one of the $K_{t-1}$ pairs of events, indicated by $p_{t-1}^{l}$  and $\ind \{D_{t-1}^l >x\}$, hold simultaneously, we have
\begin{align*}
    &D_{t}- B\wedge \left[x-D_{t-1}+B\wedge \left[\dots\left[x-D_{R+1}\right]_+\dots\right]_+\right]_+ \\
& \qquad =  D_{t}-B\wedge [x-D_{t-1}^l]_+
=  D_t - B\wedge 0
= D_t^1.
\end{align*}
The first equality is due to $p_{t-1}^{l}=1$ and the induction hypothesis. The second equality is due to $\ind \{D_{t-1}^l >x\}=1$ and the last equality is the definition of $D_t^1$.

Similarly, if the event indicated by $p_t^{K_t}$ holds, , \ie, one of the $K_{t-1}$ pairs of events, indicated by $p_{t-1}^{l}$  and $\ind \{D_{t-1}^l \le x-B\}$, hold simultaneously, we have
\begin{align*}
    &D_{t}- B\wedge \left[x-D_{t-1}+B\wedge \left[\dots\left[x-D_{R+1}\right]_+\dots\right]_+\right]_+ \\
&\qquad =  D_{t}-B\wedge [x-D_{t-1}^l]_+
=  D_t - B
= D_t^{K_t} .
\end{align*}

If the event indicated by $p_t^{k}$ holds, $k \in \Kcal_{t}\backslash \{1,K_t\}$,
\begin{align*}
    &D_{t}- B\wedge \left[x-D_{t-1}+B\wedge \left[\dots\left[x-D_{R+1}\right]_+\dots\right]_+\right]_+ \\
&\qquad =   D_{t}-B\wedge [x-D_{t-1}^{k-1}]_+
=  D_t - [x-D_{t-1}^{k-1}]
= D_t^k.
\end{align*}

As a consequence, Equation~\eqref{eq:DtkEquivalence} holds for any $t\in \Tcal$, which completes the proof.

Notice the second equality in Equation~\eqref{eq:JhatDtkForm} is another form for the same result, as the set of inequality denoted by $p_t^k$ can always be expressed as the vector inequalities. For general deficit process, this form is not useful in term of computation (and therefore we don't derive further the expression for upper and lower bounds involved). But in the setup of independent error, this forms gives computational efficient way of evaluate the thresholds as explained in Section~\ref{sec:gaussianDispatch}.
\end{IEEEproof}

\item The constrained subgradient for the cost-to-go function at stage $R+1$ given the information at a dispatch stage $r$ is
\begin{equation*}
    \grad \hat{J}_{R+1}(x_{R+1}) =-\frac{c}{T}\expec  \left\{\sum_{t\in \Tcal} \sum_{k \in \Kcal_t} \left[k \wedge(K_t+1-k)\right] p_t^k \ind \{D_t^k >x\} \middle|Y_r \right\}.
\end{equation*}
\begin{IEEEproof}
Using the explicit definition of $D_t^k$, we notice
\[
    \frac{d D_{t}^k}{d x} = \begin{cases}
    - (k-1) & \mbox{ if }1 \le k \le t-R,\\
    - (K_t - k) & \mbox{ if }t-R < k \le K_t.
    \end{cases}
\]
Consequently,
\[
    \frac{d [D_t^k-x]}{d x} =\begin{cases}
    - k &\mbox{ if } 1 \le k \le t-R,\\
    - (K_t+1 - k) & \mbox{ if }t-R < k \le K_t.
    \end{cases}
\]
or more concisely
\[
    \frac{d [D_t^k-x]}{d x} = - [k\wedge (K_t+1 - k)].
\]
Invoking the chain rule and Leibniz's rule for differentiation under the integral sign finishes the proof.
\end{IEEEproof}
\end{itemize}
\bigskip

\noindent{\bf Proof of Proposition~\ref{thm:stateSpaceDecompIndependent}}\\
By induction over the levels of the lattice. The base case $p_{R+1}^1=1$ holds trivially. Suppose the expressions for $p_{t}^k$, $\overleftarrow{p_t^k}$, $\overline{p_t^k}$ and $\overrightarrow{p_t^k}$ hold for all $t<l$ and $k\in \Kcal_t$. Consider the corresponding probabilities in level $l$. The correctness of Equation~\eqref{eq:ptkIndependent} follows from the definition of $\overleftarrow{p_{l-1}^k}$, $\overline{p_{l-1}^k}$ and $\overrightarrow{p_{l-1}^k}$. Referring to the algebraic recombinant lattice, we notice, for $k = 1$ (or $k=K_l$)
\begin{eqnarray*}
\overleftarrow{p_l^k} &=& \prob (\mbox{visit node $(l,k)$ and then move left})\\
 & = & \prob (\mbox{visit node $(l,k)$})\prob (\mbox{move left from node $(l,k)$})\\
 & = & p_{l}^{k} \prob (\epsilon_l^k >  \overline{\epsilon_l^k}).
\end{eqnarray*}
Here the last equality follows from the fact that $\epsilon_l^1=\epsilon_l^{K_l}=\epsilon_{l}$ is independent from all the past errors. Further, this result agrees with the form in Equation~\eqref{eq:ptkLeft} because $h=0$ for $k = 1$ (or $k=K_l$). Now consider $k\in \Kcal_l\backslash \{1,K_l\}$. By Equation~\eqref{eq:etk}, $\epsilon_l^k$ is independent with $\epsilon_j$ for $j<l-h$, where $h = (k-1)\wedge (K_l-k)$. Thus similar to above, we can break the joint probability into product by the independence:
\begin{eqnarray*}
\overleftarrow{p_l^k} &=& \prob (\mbox{visit node $(l,k)$ and then move left})\\
 & = & \prob (\mbox{visit node $(l-h,k-h)$}) \cdot\\
    && \prob (\mbox{starting from node $(l-h,k-h)$, visit node $(l,k)$ and then move left from it})\\
 & = & p_{l-h}^{k-h}\prob \left(\underline{\Ebold_{l-1}^{k-1}}< \Ebold_{l-1}^{k-1}\le \overline{\Ebold_{l-1}^{k-1}}, \,\, \epsilon_l^k>   \overline{\epsilon_l^k}\right).
\end{eqnarray*}
The second term in the last line follows from the observation that $k-h=1$ when $k\le (K_t+1)/2$ and $k-1 = K_{t-h}$ otherwise. That is, on the lattice, node $(l-h,k-h)$ is always a "boundary node" that is corresponding to the state either the storage is full or empty. Further, starting from such a node, the only path to node $(l,k)$ is by moving to the middle child recursively $h$ times. The exact same reasoning holds for $\overline{p_l^k}$ and $\overrightarrow{p_l^k}$ by replacing the last inequality on $\epsilon_{l}^k$ correspondingly. Thus we have proved Equations~\eqref{eq:ptkIndependent},~\eqref{eq:ptkLeft},~\eqref{eq:ptkMiddle} and~\eqref{eq:ptkRight} hold inductively. Note all the bounds in the inequalities involved are due to Corollary~\ref{thm:stateSpaceDecompGen}, with the predicted deficit term (See Equation~\eqref{eq:DtkIndependent}) plugged in.
\bigskip

\noindent{\bf Proof of Lemma~\ref{thm:offline}}\\
For the terminal cost-to-go, by Equation~\eqref{eq:JhatDtkForm},
\begin{align*}
\displaybreak
\hat{J}_{R+1}(x_{R+1})=& c \expec \left[\sum_{t\in \Tcal} \sum_{k \in \Kcal_t}[D_t^k - x]_+p_t^k\middle|Y_r\right]
= c \expec \left[\sum_{t\in \Tcal} \sum_{k \in \Kcal_t}[\hat{D}_t^k+\epsilon_t^k - x]_+ p_t^k\middle|Y_r\right]  \\
=& c  \sum_{t\in \Tcal} \sum_{k \in \Kcal_t}\bigg\{[\hat{D}_t^k - x]\overleftarrow{p_t^k} \\
 &+ \expec\left[\epsilon_t^k \middle|Y_r,\underline{\Ebold_{t-1}^{k-1}}< \Ebold_{t-1}^{k-1}\le \overline{\Ebold_{t-1}^{k-1}}, \,\, \epsilon_t^k>   \overline{\epsilon_t^k}\right]p_{t-h}^{k-h} \prob \left(\underline{\Ebold_{t-1}^{k-1}}< \Ebold_{t-1}^{k-1}\le \overline{\Ebold_{t-1}^{k-1}}, \,\, \epsilon_t^k>   \overline{\epsilon_t^k}\right)\bigg\}\\
= &c  \sum_{t\in \Tcal} \sum_{k \in \Kcal_t}\left\{\left[\hat{D}_t^k + \expec\left[\epsilon_t^k \middle|Y_r,\underline{\Ebold_{t-1}^{k-1}}< \Ebold_{t-1}^{k-1}\le \overline{\Ebold_{t-1}^{k-1}}, \,\, \epsilon_t^k>   \overline{\epsilon_t^k}\right]- x\right]\overleftarrow{p_t^k}\right\}.
\end{align*}
For the subgradient, notice $\hat{J}_{R+1}(x_{R+1})$ is a continuous function (see~\eqref{eq:JhatContReveal}). Given the Gaussian prediction error, it follows all the terms due to differentiating the integrating limits cancel. The remaining terms are given in~\eqref{eq:subgradientIndependent}.

\bibliographystyle{IEEEtran}
\bibliography{RLD}
\end{document}

%% file: defs.tex
\newcommand{\BEAS}{\begin{eqnarray*}}
\newcommand{\EEAS}{\end{eqnarray*}}
\newcommand{\BEQ}{\begin{equation}}
\newcommand{\EEQ}{\end{equation}}
\newcommand{\BIT}{\begin{itemize}}
\newcommand{\EIT}{\end{itemize}}

\newcommand{\eg}{{\it e.g.}}
\newcommand{\ie}{{\it i.e.}}

\newcommand{\ind}{\mathbbm{1}}

\newcommand{\ycal}{Y}

\newcommand{\real}{\mathbb{R}}

\newcommand{\prob}{\mathbb{P}}
\newcommand{\expec}{\mathbb{E}}

\newcommand{\grad}[1]{\nabla#1} 

\newtheorem{theorem}{Theorem}[section]
\newtheorem{lemma}[theorem]{Lemma}
\newtheorem{proposition}[theorem]{Proposition}
\newtheorem{corollary}[theorem]{Corollary}

\newtheorem{remark}[theorem]{Remark}

\usepackage{nomencl}

\makenomenclature

\usepackage{makeidx}
\makeindex

\usepackage{tikz}
\usetikzlibrary{shapes,backgrounds,calc}

\makeatletter
\tikzset{circle split part fill/.style  args={#1,#2}{%
 alias=tmp@name, 
  postaction={%
    insert path={
     \pgfextra{%
     \pgfpointdiff{\pgfpointanchor{\pgf@node@name}{center}}%
                  {\pgfpointanchor{\pgf@node@name}{east}}%
     \pgfmathsetmacro\insiderad{\pgf@x}
      \fill[#1] (\pgf@node@name.base) ([xshift=-\pgflinewidth]\pgf@node@name.east) arc
                          (0:180:\insiderad-\pgflinewidth)--cycle;
      \fill[#2] (\pgf@node@name.base) ([xshift=\pgflinewidth]\pgf@node@name.west)  arc
                           (180:360:\insiderad-\pgflinewidth)--cycle;            
         }}}}}
 \makeatother

\def\Dh{{\hat D}}

\def\D{\Delta}
\def\e{\epsilon}
\def\s{\sigma}

\usepackage{algorithmic}
\usepackage{graphicx, subfig}
\usepackage{arydshln} 